\documentclass[11pt,a4paper]{article}
\usepackage[numbers,sort&compress]{natbib}
\usepackage[colorlinks,
            linkcolor=blue,
            anchorcolor=green,
            citecolor=magenta
            ]{hyperref}
\usepackage{mathrsfs,amssymb,amsfonts}
\usepackage{graphicx}
\usepackage{amsmath,amssymb,latexsym,color}
\usepackage[mathscr]{eucal}
\usepackage[numbers]{natbib}
\usepackage{bm}
\usepackage[top=1in, bottom=1in, left=1.25in, right=1.25in]{geometry}
\usepackage{indentfirst}
\usepackage{caption}
\usepackage{subfigure}
\usepackage{setspace}%行间距
\usepackage{booktabs}            %插表格用的宏包
\usepackage{diagbox}             %插表格用的宏包
\usepackage{multirow}            %插多行表格用的宏包
\makeatletter

\newcommand{\Rmnum}[1]{\expandafter\@slowromancap\romannumeral #1@}
\makeatother
\newtheorem{theorem}{Theorem}[section]

\newtheorem{lemma}[theorem]{Lemma}
\newtheorem{remark}[theorem]{Remark}

\newtheorem{example}{Example}

\newcommand{\qed}{\hfill\Box\medskip}
\usepackage{CJK}
\usepackage{geometry}
\begin{document}
\begin{CJK*}{GBK}{song}
\begin{spacing}{1.5}
\renewcommand{\abovewithdelims}[2]{
\genfrac{[}{]}{0pt}{}{#1}{#2}}
%%%%%%%%%%%%%%%%%%%%%%%%%%%%%%%%%%%%%%%%%%%%%%%%%%%%%%%%%%%%%%%%%%%%%%%%%%%%%%%%%%%%%%%%
%%%%%%%%%%%%%%%%%%%%%%%%%%%%%%%%%%%%%%%%%%%%%%%%%%%%%%%%%%%%%%%%%%%%%%%%%%%%%%%%%%%%%%%%

\title{\bf Edge-fault-tolerance about the SM-$\lambda$ property of hypercube-like networks %\footnote{This research is supported partially by the General Project of Hunan Provincial Education Department of China (19C1742), Hu Xiang Gao Ceng Ci Ren Cai Ju Jiao Gong Cheng-Chuang Xin Ren Cai (No. 2019RS1057), the Project of Scientific Research Fund of Hunan Provincial Science and Technology Department (No. 2018WK4006).}
%\footnote{For the title, try not to use more than 3
%lines. Typeset the title in 10 pt Roman, boldface with the first letter of
%important words capitalized.}
}

%\title{\bf Edge-fault-tolerant strong Menger edge connectivity of hypercubes and folded hypercubes\\
%\footnotesize{ This research is supported by the National Natural Science Foundation of China (11571044, 61373021),  the Fundamental Research Funds for the Central University of China.} }

\author{Dong Liu \quad   Pingshan Li
 \footnote{Corresponding author. \newline {\em E-mail address:} 201921001142@smail.xtu.edu.cn(D.Liu); lips@xtu.edu.cn(P.Li); zhangbicheng@.xtu.edu.cn(B.Zhang).} \quad Bicheng Zhang \\
 {\footnotesize \em Key Laboratory of Intelligent computing$\And$Information processing of Ministry of Education}\\
{\footnotesize   \em
School of  Mathematical and Computational Sicence, Xiangtan University, Xiangtan, Hunan 41105, PR China }}
 \date{}
 \date{}
 \maketitle

\begin{abstract}

The edge-fault-tolerance of networks is of great significance to the design and maintenance of networks. For any pair of vertices $u$ and $v$ of the connected graph $G$, if they are connected by $\min \{ \deg_G(u),\deg_G(v)\}$  edge-disjoint paths, then $G$ is strong Menger edge connected (SM-$\lambda$ for short).
%A connected graph $G$ is strong Menger edge connected (SM-$\lambda$ for short) if there are $\min \{ \deg_G(u),\deg_G(v)\}$  edge-disjoint paths between any pair of vertices $u$ and $v$ of  $G$.
The conditional edge-fault-tolerance about the SM-$ \lambda$ property of $G$, written $sm_\lambda^r(G)$, is the maximum value of $m$  such that $G-F$ is still SM-$\lambda$ for any edge subset $F$ with $|F|\leq m$ and $\delta(G-F)\geq r$, where $\delta(G-F)$ is the minimum degree of $G-F$. %is the maximum integer $m$  for which $G-F$ is still SM-$\lambda$ for any edge subset $F$ with $|F|\leq m$ and $\delta(G-F)\geq r$, where $\delta(G-F)$ is the minimum degree of $G-F$.
 Previously, most of the exact value for $sm_\lambda^r(G)$ is aimed at some well-known networks when $r\leq 2$, and a few of the lower bounds on some well-known networks for $r\geq 3$.
In this paper, we  firstly determine the exact value of $sm_\lambda^r(G)$ on class of hypercube-like networks  (HL-networks for short, including hypercubes, twisted cubes, crossed cubes etc.) for a general $r$, that is,  $sm_\lambda^r(G_n)=2^r(n-r)-n$ for every $G_n\in HL_n$, where $n\geq 3$ and $1\leq r \leq n-2$.


\medskip
%\noindent {\em 2010 MSC:} 05C25; 05C15.
\end{abstract}

\section{Introduction}
No matter what kind of interconnection network, there are growing concern over the fault-tolerance and reliability. The edge connectivity is one of the indexes  extensively used to measure the fault-tolerance of networks. In graph-theoretic terms, the edge connectivity of  the graph $G$, written  $\lambda(G)$, is the minimum number of  edges whose  removal disconnected the graph $G$. In contrast to edge connectivity as a global concept, Menger theorem provides a local perspective, which uses the number of edge-disjoint paths between two vertices to express the edge connectivity between them.  The edge-disjoint paths of networks can not only  improve the efficiency and security of network transmission, but more importantly, it can provide fault-free paths when the network fails. Thus, the multiplicity of edge-disjoint paths is also a suitable measure for fault-tolerance of networks.  Usually,  it's better if there are more edge-disjoint paths between two vertices in a graph.
%the more edge-disjoint paths, the better  fault-tolerance of  networks.

However, the number of edge-disjoint paths between two vertices is limited by the degree of these two vertices.  In other words, for any pair of vertices $u$ and $v$ of a graph $G$,  there are at most $\min\{ \deg_G(u),\deg_G(v)\}$  edge-disjoint paths between $u$ and $v$, where $\deg_G(u)$ and $\deg_G(v)$ are degrees of vertices $u$ and $v$, respectively. When designing a network, it is desirable to maximize the number of edge-disjoint paths between two vertices. %It is desirable to have as many edge-disjoint paths as possible between two vertices.
 Thus, Fricke et al.{\rm\cite{fricke}}  begun to take an interest in the strongly Menger edge connected (SM-$\lambda$ for short) graph in which there are $\min\{ \deg_G(u),\deg_G(v)\}$ edge-disjoint paths for every pair of vertices $u$ and $v$ of the graph.
%where any pair of vertices $u$ and $v$ are connected by $\min\{ \deg_G(u),\deg_G(v)\}$ edge-disjoint paths, which is called strongly Menger edge connected (SM-$\lambda$ for short) graph.
Sometimes it is called the maximally local edge connected graph. For more information about the SM-$\lambda$ graph, please refer to the papers {\rm \cite{hellwig,holtkamp,volkmann}}.

Once the scale of the network becomes larger,  network failures are inevitable. Therefore, the research on faulty networks is necessary. There are different models based on the distribution of fault edges. One of them is the random faulty model in which the fault edges can be randomly distributed without restriction. This model is applied in {\rm \cite{zhai,Li,gu1,Li1,guo}}. However, it is  rare that almost all edges incident to a vertex fail simultaneously. Thus, a conditional faulty model was proposed, under which each vertex is incident to at least two fault-free edges. This model has been applied in {\rm \cite{qchen,qiao,Li,gu1,Li1,guo}}. Later, He et al.{\rm\cite{he}} extended this model. Under the new conditional faulty model, each vertex is incident to at least $r$ fault-free edges.  This model is applied in {\rm \cite{zhai,gu}}. And by the way, the random faulty model is the new conditional faulty model with $r=0$.

There is an interesting question: what is the maximum number of faulty edges such that the graph is still SM-$\lambda$ under the  new condition faulty model? In response to this problem, He et al.{\rm\cite{he}} proposed a parameter to consider the edge-fault-tolerance about the SM-$\lambda$ property of a graph. The conditional edge-fault-tolerance about the SM-$ \lambda$ property of $G$, written $sm_\lambda^r(G)$, is the maximum value of $m$  such that $G-F$ is still SM-$\lambda$ for any edge subset $F$ with $|F|\leq m$ and $\delta(G-F)\geq r$, where $\delta(G-F)$ is the minimum degree of $G-F$. Recently,  the edge-fault-tolerance about the SM-$\lambda$ property of some famous networks has been studied, %the parameter $sm_\lambda^r(G)$  have been studied for several famous networks,
such as $n$-dimensional hypercube $Q_n$, $n$-dimensional folded hypercube $FQ_n$, $n$-dimensional balanced hypercube $BH_n$, augmented $k$-ary $n$-cube $AQ_{n,k}$, $n$-dimensional bubble-sort star graph $BS_n$, $(n,k)$-star graph $S_{n,k}$, data center network $ D_{n,k}$, $n$-dimensional augmented cube $AQ_n$, $n$-dimensional hypercube-like network (HL-network for short) $G_n$.  For more details, please refer to the Table \ref{tab1}.

So far, most of the papers only gave the exact value of $sm_\lambda^0(G)$ and $sm_\lambda^2(G)$ for some specific networks. For $r\geq 3$, only the lower bound of $sm_\lambda^r(G)$ of some networks was given except $sm_\lambda^3(AQ_n)=6n-15$ \rm{\cite{zhai}}. For $n$-dimensional HL-network $G_n$, Li and Xu \cite{Li1} determined $sm_\lambda^0(G_n)$ and $sm_\lambda^2(G_n)$. In this paper, we show that $sm_\lambda^r(G_n)=2^r(n-r)-n$ for $1\leq r \leq n-2$, where $n\geq 3$.

The rest of this paper is organized as follows. In Section $2$, we  introduce the definition and some properties of HL-networks. In Section $3$, we determine the exact values of  $sm_\lambda^r(G_n)$ for the  $n$-dimensional HL-network $G_n$.
\begin{table}[htbp]
\renewcommand\arraystretch{1.2}
\renewcommand\tabcolsep{8.0pt}
%h：hear，t：top，b：bottom，p：page，下一页。
\centering
\caption{$ sm_\lambda^r(G)$ of some famous  networks}
\label{tab1}
    \begin{tabular}{cccc}
    %几列就写几个c，表示内容居中写，哪里需要分割线就在哪里加“|”
        \toprule
         %最上面的横线
        Graph & $ \delta(G-F)\geq r$ & $  sm_\lambda^r(G)$ & ref\\
        %填写表格内容，建议内容多的时候用{} & {} & {}
        \midrule
        %行间需要分割线就加这条语句，不需要就不加
        %只要是行间的内容，填充完之后都要加上“\\”
     & $r=2$& $\geq2n-4$ ($n\geq4$) & \rm{\cite{qiao}} \\
       Hypercube &$r=3$ & $\geq 3n-9$ $(n\geq7)$&\rm{\cite{he}} \\
        $Q_n$&$r=4$ &$\geq4n-14 $ $(n\geq8)$&\rm{\cite{he}}\\
        &$ r=5$&$\geq5n-20$ $(n\geq8)$ &\rm{\cite{he}}\\
        \hline
            Folded hypercube&\multirow{2}{*}{$r=2$} &$\geq2n-2$ ($n\geq5$)& \rm{\cite{qiao}}\\
       $FQ_n$& &$=3n-5$ ($n\geq5$)& \rm{\cite{qchen}} \\
         %\midrule
         \hline
          Balanced hypercube & $r=0$ & $=2n-2$&  \multirow{2}{*}{ \rm{\cite{Li}}}\\
         $BH_n$ &$r=2$ &$=6n-8$ ($n\geq2$) &\\
          %\midrule
          \hline
         Augmented $k$-ary &$r=0$ &$=4n-4$ ($n\geq2,k\geq 3$)& \multirow{2}{*}{\rm{\cite{gu1}}}\\
         $n$-cube $AQ_{n,k}$&$r=2$&$=8n-10$ ($n\geq2, k\geq 3$)&  \\
         \hline
        Bubble-sort star graph &$r=0$ &  $=2n-5$ $(n\geq3)$& \multirow{2}{*}{ \rm{\cite{guo}}} \\
         $BS_n$&$r=2$ &$=6n-17$ $(n\geq4)$ &\\
         \hline
     Data center network  $D_{n,k}$  &$1\leq r\leq 3$ & $\geq rn+rk-r^2-2r\ (n\geq3,k\geq4)$ &\multirow{2}{*}{\rm{\cite{gu}}} \\
        $(n,k)$-star graph $S_{n,k}$  &$1\leq r\leq 3$  & $\geq rn-r^2-2r\ (n\geq r+3,k\geq3)$ &  \\
         \hline
          \multirow{2}{*}{Augmented cube } &$r=0$ & $=2n-3$ ($n\geq4$)&\multirow{3}{*}{\rm{\cite{zhai}}} \\
         \multirow{2}{*}{$AQ_n$} &$r=3$ &$=6n-15$ ($n\geq4$) &\\
        &$r\geq3$  & $\geq2nr-4r(r+1)+8$ ($n\geq8$)& \\
        \hline
         \multirow{2}{*}{Hypercube-like network}  &$r=0$ & $=n-2$ ($n\geq4$)& \rm{\cite{Li1}}\\
        \multirow{2}{*}{$G_n$} &$r=2$ &$=3n-8$ ($n\geq3$) &\rm{\cite{Li1}}\\
         &$1\leq r\leq n-2 $ &$=2^r(n-r)-n \ (n\geq3)$ & our work\\
        \bottomrule
        %最下面的横线
\end{tabular}
\end{table}

\section{Some properties of  HL-networks}
As one of  popular topologies for the design of the computer system, the hypercube $Q_n$  has many excellent properties such as good connectivity, recursive scalability and symmetry. However, the diameter of the hypercube is relatively large. Efforts to improve this shortcoming have led to the emergence of many hypercube variants,  such as twisted cubes \rm{\cite{Hibers}},  crossed cubes \rm{\cite{Efe}} and M$\rm{\ddot{o}}$bius cubes  \rm{\cite{Cull}}. They not only have a smaller diameter than the hypercube, but also retain some of the good properties of the hypercube. To conduct integrated research on these variants, hypercube-like networks (HL-networks for short) were proposed by Vaidya et al.\rm{\cite{Vaidya}}. In \rm{\cite{Fan}}, HL-networks is also called bijective connection networks. There are many research on HL-networks, for example {\rm \cite{Chang,Chang2,Li,Park,Zhou,shi}}).

\begin{figure}[htbp]
\centering	
	\subfigure[]
	{
		\begin{minipage}{6cm}
			\centering
			\includegraphics[scale=0.4]{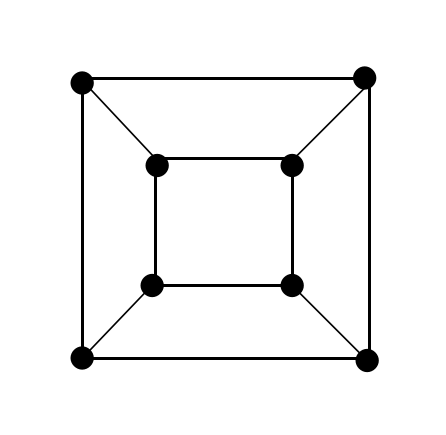}
		\end{minipage}
	}
	\subfigure[]
	{
		\begin{minipage}{6cm}
			\centering
			\includegraphics[scale=0.4]{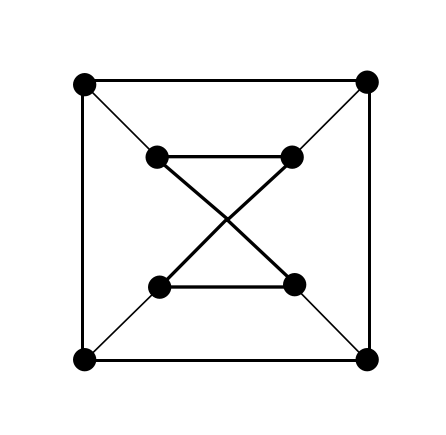}
		\end{minipage}
	}
	\caption{Two 3-dimensional HL-networks: (a) $Q_3$\ \  (b) $CQ_3$}
    \label{f1}
\end{figure}

 Let $HL_n$ be the set of all $n$-dimensional HL-networks and  $G_n$ be an any element of $HL_n$. The recursive definition of HL-networks is as follows:

\begin{enumerate}
\item $HL_0=\{K_1\}$, where $K_1$ is an isolate vertex;
\item $HL_{n}=\{G_{n-1}\oplus G'_{n-1}|G_{n-1},G'_{n-1}\in HL_{n-1}\}$,
 \end{enumerate}
where the symbol $``\oplus"$ represents a perfect matching operation between two $(n-1)$-dimensional subcubes $G_{n-1}$ and $G'_{n-1}$.
 Clearly, $HL_1=\{K_2\}$, $HL_2=\left\{ C_4\right\}$, and $HL_3=\left\{ Q_3,\  CQ_3\right\}$, where $C_4$ is a cycle of length 4, and $Q_3$ and $CQ_3$ are shown in Figure \ref{f1}.  By the definition, every $G_n$ is an $n$ regular graph with $2^n$ vertices and $n2^{n-1}$ edges.

Next, we will introduce some properties of HL-networks.

Let $g$ be an integer. Then $g$ can be written as the sum of the exponents of 2, that is, $g=\begin{matrix} \sum_{i=0}^s 2^{t_i}\end{matrix} $, where $t_0 = \left\lfloor \log_2 g\right \rfloor$, $t_i = \left\lfloor \log_2(g- \begin{matrix} \sum_{r=0}^{i-1} 2^{t_r}\end{matrix})\right\rfloor $ for $i \ge 1$. For any $G_n\in HL_n$, we use $e_g$  to denote the maximum number of edges of the subgraph (of $G_n$) induced by $g$ vertices.
For any $X\in V(G_n)$, let $E_X$ be a set of edges where each edge has exactly one endpoint in $X$. For convenience, we define
\begin{center}
$f(g)=ng-2e_g$.
\end{center}
\begin{lemma}\rm{\cite{liu}}\label{le2.1}
Given any $G_n\in HL_n$, if $g=\sum_{i=0}^s2^{t_i}$, then $e_g=\sum_{i=0}^st_i2^{t_i-1}+\sum_{i=0}^{s}i2^{t_i} $.
\end{lemma}
\begin{lemma}\rm{\cite{liu}}\label{le2.2}
For any $X\in V(G_n)$ with $|X|=g$ , $|E_X|\geq ng-2e_g$.
\end{lemma}
\begin{lemma}\rm{\cite{Fan}}\label{le2.3}
For any $G_n\in HL_n$, $\lambda(G_n)=n$ for $n\geq 1$.
\end{lemma}
\begin{lemma}\label{le2.4}
If $1\leq g\leq 2^n-1$, then $f(g)\geq n$.
\end{lemma}
\noindent{\bf Proof: }Take a connected subgraph $H$ from $G_n$ such that $|V(H)|=g$ and $|E(H)|=e_g$. Let $S=E_{V(H)}$, then $|S|=n|V(H)|-2|E(H)|= ng-2e_g=f(g)$. Because $1\leq g \leq2^n-1$, $S$ is an edge cut of $G_n$. By Lemma \ref{le2.3}, $f(g)\geq n$.$\qed$
\begin{lemma}\rm{\cite{liu}}\label{le2.5}
If $g\leq 2^{n-2}$, then $(n-2)g-2e_g\geq 0$.
\end{lemma}
\begin{lemma}\label{le2.6}
Let $r$ be an integer. If $2^r\leq g\leq 2^{n-1}$, then $f(g)\geq f(2^r)$,  where $0\leq r\leq n-1$.
\end{lemma}
\noindent{\bf Proof: } We firstly claim that
\begin{center}
 $f(g)\geq f(2^{r+i})$ when $2^{r+i}\leq g < 2^{r+i+1}$ ($0\leq i\leq n-r-2$).
\end{center}

Let $g=\sum^s_{i=0}2^{t_i}$, then $t_0=r+i$ and $t_1\leq r+i-1\leq r+n-r-2-1=n-3$. Thus,
\begin{center}$\begin{array}{rl}
f(g)-f(2^{r+i})&=[n\sum^s_{i=0}2^{t_i}-2(\sum^s_{i=0}t_i2^{t_i-1}+\sum^s_{i=0} i2^{t_i})]-[n2^{r+i}-2(r+i)2^{r+i-1}]\vspace{1.0ex}\\
%&=n(2^{r+i}+\sum^s_{i=1}2^{t_i})-2[(r+i)2^{r+i-1}+\sum^s_{i=1}t_i2^{t_i-1}+\sum^s_{i=1} i2^{t_i}]-[n2^{r+i}-2(r+i)2^{r+i-1}]\vspace{1.0ex}\\
&=n\sum^s_{i=1}2^{t_i}-2(\sum^s_{i=1}t_i2^{t_i-1}+\sum^s_{i=1} i2^{t_i})\vspace{1.0ex}\\
&=(n-2)\sum^s_{i=1}2^{t_i}-2[\sum^s_{i=1}t_i2^{t_i-1}+\sum^s_{i=1} (i-1)2^{t_i}]\vspace{1.0ex}\\
&=(n-2)\sum^s_{i=1}2^{t_i}-2e_{(\sum^s_{i=1}2^{t_i})}.
\end{array}$
\end{center}
Note that $t_1\leq n-3$, then $\sum^s_{i=1}2^{t_i}< 2^{n-2}$. By Lemma \ref{le2.5},  $(n-2)\sum^s_{i=1}2^{t_i}-2e_{(\sum^s_{i=1}2^{t_i})}\geq 0$. Thus, $f(g)\geq f(2^{r+i})$.

Next, we prove that
\begin{center}
$f(2^{r+i+1})\geq f(2^{r+i})$ for $0\leq i\leq n-r-2$.
\end{center}
Note that $e_{2^{r+i+1}}=(r+i+1)2^{r+i}$, then we have that
\begin{center}$\begin{array}{rl}
f(2^{r+i+1})&=n2^{r+i+1}-2(r+i+1)2^{r+i}\vspace{1.0ex}\\
&=2^{r+i}[n-(r+i)]+2^{r+i}(n-r-i-2)\vspace{1.0ex}\\
&=n2^{r+i}-2(r+i)2^{r+i-1}+2^{r+i}(n-r-i-2)\vspace{1.0ex}\\
&=f(2^{r+i})+2^{r+i}(n-r-i-2).\vspace{1.0ex}\\
\end{array}$
\end{center}
Thus,  $f(2^{r+i+1})\geq f(2^{r+i})$ when $0\leq i\leq n-r-2$.

We complete the proof as follows.
\begin{itemize}
\item If $2^r\leq g < 2^{n-1}$, then there exists an integer $j$ $(0\leq j\leq n-r-2)$ such that $2^{r+j}\leq g< 2^{r+j+1}$. Thus,
\begin{center}
$f(g)\geq f(2^{r+j})\geq f(2^{r+j-1})\geq...\geq f(2^r)$.
\end{center}
\item If $g=2^{n-1}$, then $f(g)=f(2^{n-1})\geq f(2^{n-2})\geq...\geq f(2^r)$.$\qed$
\end{itemize}
\begin{lemma}\label{le2.7}
 Given any $G_n\in HL_n$ %let $S$ be an arbitrary subset of $E(G_n)$ with $|S|\leq 2^r(n-r)-1$,
 and edge subset $S\in E(G_n)$ with $|S|\leq 2^r(n-r)-1$, then there exists a component $C$ in $G_n-S$ such that $|V(C)|\geq 2^n-(2^r-1)$, where $n\geq3$ and $0\leq r\leq n-2$.
\end{lemma}
\noindent{\bf Proof: }Let $C_1,C_2,...,C_k$ be components of $G_n-S$ such that $|V(C_1)|\leq|V(C_2)|\leq...\leq|V(C_k)|$, then we need to show that $|S|\geq 2^r(n-r)$ if $|V(C_k)|\leq2^n-2^r$.

\noindent{\bf Case 1: }Suppose that $|V(C_k)|<2^r$.

In this case, $|V(C_i)|<2^r$ for every $i\in\{1,2,...,k\}$. Note that $\sum_{i=0}^k|V(C_i)|=2^n$, then there exists an integer
$j$ such that $\sum_{i=1}^j |V(C_i)|< 2^r$ but $2^{r}\leq\sum_{i=1}^{j+1}|V(C_i)| <2^{r+1}$. Let $X_1=\bigcup_{i=1}^{j+1}V(C_i)$. Clearly, $2^r\leq |X_1|<2^{r+1}$. Note that there are more than two components  in $X_1$.  Thus,
\begin{center}
$|S|>|E_{X_1}|\geq n|X_1|-2e_{|X_1|}\geq f(2^r)=2^r(n-r)$.
\end{center}
The last inequality holds by Lemma \ref{le2.6}.

\noindent{\bf Case 2: }Suppose that $2^r\leq |V(C_k)|< 2^{n-1}$.

 By Lemma \ref{le2.6}, we have that
\begin{center}
$|S|>|E_{V(C_k)}|\geq n|V(C_k)|-2e_{|V(C_k)|}\geq f(2^r)=2^r(n-r)$.
\end{center}

\noindent{\bf Case 3: }Suppose that $2^{n-1}\leq |V(C_k)|\leq 2^n-2^r$.

 In this case, $2^r\leq \sum_{i=1}^{k-1}|V(C_i)|\leq 2^{n-1}$. Let $X_3=\bigcup_{i=1}^{k-1}V(C_i)$. Clearly, $2^r\leq |X_3|\leq2^{n-1}$. Thus,
\begin{center}
$|S|\geq |E_{X_3}|\geq n|X_3|-2e_{|X_3|}\geq f(2^r)=2^r(n-r)$.
\end{center}
The last inequality holds by Lemma \ref{le2.6}.$\qed$

\section{\bf $ sm_\lambda^r(G)$ for HL-networks }
For $n$-dimensional HL-network $G_n$, Li and Xu \cite{Li1} determined $sm_\lambda^0(G_n)$ and $sm_\lambda^2(G_n)$ in \cite{Li1}. In this section, we determine the exact value of  $sm_\lambda^r(G_n)$ for a general integer $r$.

Let $H_1$ and $H_2$ be  disjoint subgraphs of the graph $G$ and $u\in V(G)$. The set of all edges between $H_1$ and $H_2$ is written as $E(H_1, H_2)$.  The set of all vertices adjacent to $u$ is written as $N_G(u)$.
%The set of all edges between $G_1$ and $G_2$ is denoted by $E(G_1, G_2)$.  The set of all neighbors of $u$ is denoted by $N_G(u)$. Similarly, the set of all neighbors of all vertices in $G_1$ is denoted by $N_G[V(G_1)]$.
For any pair of vertices $u$ and $v$ of the graph $G$, the $u$,$v$-edge cut is the edge set whose deletion makes $u$ and $v$  disconnected.

\begin{theorem}\rm{\cite{menger}}\label{le3.1}
Let $u$ and $v$ be a pair of vertices in the graph $G$. The minimum size of an $u$, $v$-edge cut equals the maximum number of edge-disjoint $u$, $v$-paths.
\end{theorem}
\begin{theorem}
 For any $G_n\in HL_n$, $sm_\lambda^r(G_n)=2^r(n-r)-n$, where $n\geq 3$ and $1\leq r\leq n-2$.
\end{theorem}
\noindent{\bf Proof: }First,we prove that $sm_\lambda^r(G_n)\geq 2^r(n-r)-n$.  Suppose that  $F$ is an %faulty
edge set of $G_n$ such that $|F|\leq 2^r(n-r)-n$ and $\delta(G_n-F)\geq r$. Now we prove that $ G_n -F$ is still SM-$\lambda$.

We claim that $G_n-F$ is connected. Conversely, suppose that $G_n-F$ is not connected. Since $|F|\leq 2^r(n-r)-n<2^r(n-r)-1$, there is a component $C$  with at least $2^n-(2^r-1)$ vertices in $G_n-F$ by Lemma \ref{le2.7} . Let $|V(G_n-F-C)|=t$. Clearly, $1\leq t\leq 2^r-1$. Since $\delta(G_n-F)\geq r$, $|E(G_n-F-C)|\geq \frac{rt}{2}$. Note that $|E(G_n-F-C)|\leq e_t$. We have that

$$\frac{rt}{2}\leq |E(G_n-F-C)|\leq e_t. \eqno{(1)}$$
Because $1<t<2^r-1$, one has $f(t)=rt-2e_t\geq r>0$ by Lemma \ref{le2.4}. This contradicts (1).  Thus,  the claim holds.

Let $u$, $v$ be any pair of vertices of $G_n$ such that $\deg_{G_n-F}(u)\leq \deg_{G_n-F}(v)$.  By the definition and Theorem \ref{le3.1}, we need to prove that removing at most $\deg_{G_n-F}(u)-1$ edges from $G_n-F$ is not enough to make  $u$ and $v$ disconnected
%we  need to prove $u$ and $v$ are still connected if the number of edges removed is no more than $\deg_{G_n-F}(u)-1$ in $G_n-F$.

Conversely, suppose that $u$ and $v$ are disconnected after removing an edge set from $G_n-F$ whose size is $\deg_{G_n-F}(u)-1$ at most. Let $E_f$ be the minimum one. Since $\deg_{G_n-F}(u)\leq n$,  one has  $|E_f|\leq n-1$. Suppose that $S=E_f\cup F$. We have that $|S|\leq 2^r(n-r)-1$. There exists a component $T$ in $G_n-S$ and $|V(T)|\geq 2^n-(2^r-1)$ by Lemma \ref{le2.7}.
%Conversely, suppose that $u$ and $v$ are disconnected after removing a edge set from $G_n-F$.  Let $E_f$ be the minimum one, where $|E_f|\leq \deg_{G_n-F} (u)-1$. Note that $\deg_{G_n-F}(u)\leq \deg_{G_n}(u)=n$, so we have  $|E_f|\leq n-1$. Let $S=E_f\cup F$. Then, $|S|\leq 2^r(n-r)-1$. By Lemma \ref{le2.7}, there is a component $T$ in $G_n-S$ such that $|V(T)|\geq 2^n-(2^r-1)$.
Because $u$ and  $v$ are disconnected,  $u$ and $v$ are not located in the same component of $G_n-S$.  Without loss of generality, let $u$ be located in the other component $T'$ of $G_n-S$ where $|V(T')|=g$. Clearly, $1\leq g \leq 2^r-1$. Let $\deg_{T'}(u)=t$ and  $X=\{x\in V(T')\backslash \{u\}|\deg_{T'}(x)\leq r-1\}$. Then
\begin{center}
$t+\sum_{x\in X}\deg_{T'}(x)+r(g-|X|-1)\leq \sum_{x\in V(T')}\deg_{T'}(x)\leq 2e_g\leq r(g-1)$.
\end{center}
The last inequality holds by Lemma \ref{le2.4}. Thus,
\begin{center}
$t\leq r|X|-\sum_{x\in X}\deg_{T'}(x)=\sum_{x\in X}(r-\deg_{T'}(x))$.
\end{center}
Note that $\delta(G_n-F)\geq r$, then every vertex in $X$ is incident with at least $r-\deg_{T'}(x)$ edges in $E_f$. By the minimality of $E_f$, one has $E_f\cap E(G_n[X])=\emptyset$. It follows that
\begin{center}
$|E_f|\geq (\deg_{G_n-F}(u)-t)+\sum_{x\in X}(r-\deg_{T'}(x))\geq \deg_{G_n-F}(u)$.
\end{center}
It is  contradicty. Thus, $sm_\lambda^r(G_n)\geq 2^r(n-r)-n$.

Next, we  show that $sm_\lambda^r(G_n)\leq 2^r(n-r)-n$. Suppose that $G_r$ is an $r$-dimensional subcube of $G_n$. Let $u\in V(G_r)$, $N_{G_r}(u)=\{u_1, u_2,..., u_r\}$ and $e_i\in E(u_i, G_n-G_r)$. Let
\begin{center}
$F=E(G_r-\{u\}, G_n-G_r)\setminus \bigcup^{r-1}_{i=1}\{e_i\}$.
\end{center}
Note that $G_r$ is $r$ regular. Thus,
\begin{center}
$|F|=(2^r-1)(n-r)-(r-1)=2^r(n-r)-n+1$.
\end{center}

\begin{figure}[htbp]
    \centering
    \includegraphics[scale=0.5]{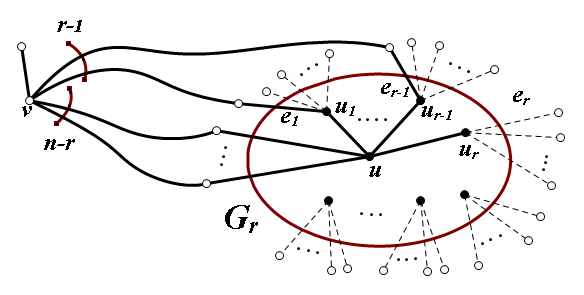}
    \caption{Illustration for Theorem 3.2}
    \label{f2}
    \end{figure}

Because $2^n-2^r(n-r+1)>1$ for $1\leq r \leq n-2$ and $n\geq3$, one can take a vertex $v\in V(G_n)\setminus \bigcup_{x\in V(G_r)}N_{G_n}(x)$. It follows that  $\deg_{G_n-F}(u)=\deg_{G_n-F}(v)=n$.  By observing Figure \ref{f2},  we can see that $u$ and $v$ are connected by $(n-1)$-edge-disjoint paths at most.
 %By observation (see Figure \ref{f2}), there are at most $(n-1)$-edge-disjoint paths between $u$ and $v$.
Then $G_n-F$ is not SM-$\lambda$. Thus, $sm_\lambda^r(G_n)\leq 2^r(n-r)-n$.$\qed$

\section{Conclusion}
In this paper, we  show that $sm_\lambda^r(G_n)=2^r(n-r)-n$ for $n\geq3$ and $1\leq r \leq n-2$, which generalizes the result of Li and Xu. More importantly, our result is the first exact value of $sm_\lambda^r(G)$ for a general integer $r$. In the future, we will consider the exact value of  $sm_\lambda^r(G)$ of other networks for a general integer $r$. In addition, we will also consider networks with vertex failures.
\end{spacing}

\end{CJK*}

\end{document}